\documentclass{l4dc2027}

\title[Learning Morse Lyapunov functions]{Learning Morse Lyapunov functions through invertible neural networks}
\usepackage{times}
\usepackage{mathrsfs}
\usepackage{enumerate}

\author{%
 \Name{Wouter Jongeneel} \Email{wouterjo@kth.se}\\
 \addr KTH Royal Institute of Technology and Digital Futures
}

\begin{document}

\maketitle

\begin{abstract}
Towards data-driven stability certificates of dynamical systems, we are witnessing renewed interest in studying and exploiting the structure of Lyapunov functions. Especially Lyapunov functions that are locally quadratic (\textit{i.e.}, Morse) gained popularity ever since it was understood that they can be written as the norm squared of a diffeomorphism. This topological representation alludes to parametrizations using invertible neural networks. In this note we discuss when and how such functions can be learned. In particular, using neural ODEs (NODEs), we characterize Sobolev approximation of these diffeomorphisms on compact sets and show how this naturally relates to stability. We also connect our results to PLNets and supply numerical experiments for illustration. 
\end{abstract}


 \section{Introduction}
After pioneering work by \cite{ref:Massera1956}, \cite{ref:kurzweil1963inversion} and \cite{ref:wilson1969smoothing}, there is now an extensive theory describing when the existence of a continuous or smooth Lyapunov function $V$ is equivalent to the existence of an asymptotically stable attractor under some dynamical system, \textit{e.g.}, see~\cite{hahn1967stability,ref:bhatia1970stability,ref:bacciotti2005liapunov,ref:fathi2019smoothing}. 

Towards a data-driven stability analysis, however, one hopes to say more beyond regularity of those Lyapunov functions, \textit{e.g.}, to work with suitable parametric function spaces when searching for a Lyapunov function, see \cite{ref:manchester2026neural} for a recent overview and \cite{ref:zkolter2019learning} for work that propelled this field of study. In this article we focus solely on, structure of the Lyapunov function, but note that a lot of work aims to study or identify some unknown dynamical system together with such a stability certificate, \textit{e.g.}, see \cite{ref:harting2026locally}. 

In a companion note {\citep{ref:jongeneel2026global}}, we discuss how the structure of Lyapunov functions helps in uncovering topological features of dynamical systems. We exploit those results below.

\subsection{Setting}
In this note we work with $C^r$-smooth dynamical systems of the form 
\begin{equation}
\label{equ:f}
    \dot{x}=f(x), \quad x\in \mathbb{R}^n,
\end{equation}
for some $r\in \mathbb{N}_{\geq 1}\cup \{\infty\}$, and such that some equilibrium point $x^*\in \mathbb{R}^n$ is \textbf{\textit{globally asymptotically stable}} (GAS). To aid the reader, we briefly recall global asymptotic stability, (\textit{e.g.}, see \cite[p. 28]{ref:bacciotti2005liapunov}). Consider~\eqref{equ:f}, and let the corresponding semiflow be denoted by $\phi$, \textit{i.e.}, $\frac{\mathrm{d}}{\mathrm{d}\tau}\phi^{\tau}(x)|_{\tau=t} = f(\phi^t(x))$ for all $(t,x)\in \mathbb{R}_{\geq 0}\times \mathbb{R}^{n}$. We will not be concerned with existence of solutions. A point $x^{*}\in \mathbb{R}^n$ is an \textit{equilibrium point} of $f$ when $f(x^{*})=0$ or equivalently, $x^*$ is a fixed point of the time-$t$ map: $\phi^t(x^{*})=x^{*}$ for any $t\geq 0$. Then, the equilibrium point $x^*$ is said to be \textit{{globally asymptotically stable}} (GAS) (under $f$) when: $(i)$ $x^*$ is \textit{Lyapunov stable}, that is, for any open neighbourhood $U_{\varepsilon}\ni x^*$ there is an open set $U_{\delta}\subseteq U_{\varepsilon}$ such that $\phi^t(U_{\delta})\subseteq U_{\varepsilon}$ $\forall t\geq 0$; and $(ii)$ $x^*$ is \textit{globally attractive}, that is, $\lim_{t\to+\infty}\phi^t(x_0)=x^*$ for all $x_0\in \mathbb{R}^n$.

It is hard to work directly with the definition of GAS. Luckily, $x^*$ being GAS is equivalent to the existence of a smooth (\textit{i.e.}, $C^{\infty}$-smooth) function $V:\mathbb{R}^n\to \mathbb{R}_{\geq 0}$ that is: $(i)$ coercive (\textit{i.e.}, all sublevel sets $\{x\in \mathbb{R}^n:V(x)\leq c\}$ are compact); $(ii)$ such that $\{x^*\}=V^{-1}(0)$; and $(iii)$ away from $x^*$, decreasing strictly under solutions of~\eqref{equ:f}, \textit{i.e.}, $L_f V(x):=\langle \nabla V(x),f(x) \rangle < 0$ for all $x\in \mathbb{R}^n\setminus \{x^*\}$, \textit{e.g.}, see \cite[Thm~2.4]{ref:bacciotti2005liapunov}. Functions $V$ satisfying properties $(i)$, $(ii)$ and $(iii)$ will be simply called smooth \textbf{\textit{Lyapunov functions}} (this instead of smooth, strict, proper Lyapunov functions with respect to equilibria). Without loss of generality, in the remainder of this note we focus on $x^*=0$. We also remark that we work on $\mathbb{R}^n$ for simplicity. 

Regarding further notation and structure, following \cite{ref:kvalheim2025differential}, let $\mathcal{S}^r_0(\mathbb{R}^n)$ be the set of $C^r$-smooth vector fields on $\mathbb{R}^n$, for some $r\in \mathbb{N}_{\geq 1}\cup \{\infty\}$, such that $x^*=0$ is GAS.

Suppose now that we have a family of Lyapunov function candidates, parametrized through $\theta\in \Theta\subseteq \mathbb{R}^d$ via some neural network: $\{ V_{\theta}\in C^{\infty}(\mathbb{R}^n;\mathbb{R}_{\geq 0}):\theta \in \Theta \}$. With the aforementioned terminology in place, the central question of this note is to understand precisely when and how these functions can be selected to be \textit{Morse}. Better yet, we discuss why this is interesting.

We recall that a function $v\in C^2(\mathbb{R}^n;\mathbb{R})$ is a \textbf{\textit{Morse function}} when all its critical points (\textit{i.e.}, $x^*\in \mathbb{R}^n$ such that $\nabla v(x^*)=0$) are non-degenerate, that is, the Hessian $\nabla^2 v(x^*)$ is non-singular. Inspired by a rich set of topological results from Morse theory (\textit{e.g.}, see \cite{ref:milnor1963morse}), one may motivate the study, \textit{globally}, of Lyapunov functions that are Morse, \textit{i.e.}, $\langle q,\nabla^2V(x^*) q\rangle > 0$ for all $q\in \mathbb{R}^n\setminus \{0\}$ and some equilibrium $x^*$. The fact that the Hessian is symmetric positive definite is written as $\nabla^2 V(x^*)\succ 0$. Such a Lyapunov function will be called a \textit{\textbf{Morse Lyapunov function}}. 

\textbf{Contribution:} In our companion note, inspired by \cite{ref:boumal2026smooth}, we studied topological ramifications of working with Lyapunov functions that are Morse \citep{ref:jongeneel2026global}. In this work, being inspired by advances in invertible neural networks (\textit{e.g.}, \cite{ref:ishikawa2023universal}) we detail effectiveness of employing invertible neural networks to learn Morse Lyapunov functions as stability certificates. In particular, we provide (i) an obstruction for the existence of Morse Lyapunov functions (Proposition~\ref{prop:obstruction}); (ii) as the obstruction is mild, we can provide theory for when NODE-based INNs yield Morse Lyapunov functions on compact sets (Theorem~\ref{thm:main}, Corollary~\ref{cor:main}); and (iii) we link these results to P\L{} Lyapunov functions, providing context for the distortion parameter $\tau$ (Remark~\ref{rem:PL}).  

\subsection{Simple Lyapunov functions do not suffice, in general}

First, when~\eqref{equ:f} is linear, $V$ can be chosen to be \textit{quadratic} (\textit{e.g.}, $x\mapsto \langle x,Px\rangle$ for $P\succ 0$), globally, \textit{e.g.}, see \cite[Thm.~18]{ref:sontag2013mathematical}. Under hyperbolicity, some arguments survive, \textit{locally}, for nonlinear systems, \textit{e.g.}, see \cite[\S 9]{ref:lefschetzlasalle1961} and \cite[\S 28]{hahn1967stability} for early results or~\cite{ref:barreau2024supervised} and \cite{ref:harting2026locally} for recent learning-oriented results.

Globally, one may hope that more can be said beyond $C^r$-smoothness of $V$, \textit{e.g.}, one may hope that any Lyapunov function can be chosen to be a polynomial. This hope was refuted by \cite{ref:ahmadi2011globally}. In a similar vein, there are vector fields that render the origin GAS, but provably fail to admit a \textit{convex} Lyapunov function \cite[Ex. III.3]{ref:JongeneelSchwan2024TAC}. Moreover, one can find examples where a smooth Lyapunov function exists, but not one that is locally real-analytic, see \cite[Prop. 5.2 (b)]{ref:bacciotti2005liapunov} and \cite{ref:ahmadi2018globallyP}.

One may continue and hope that Lyapunov functions can always be chosen such that they are Morse, for instance. With the proposition below we illustrate that in general, not much more can be said beyond simple $C^r$-regularity. The example to have in mind is a planar ODE of the form: 
    \begin{equation}
        \label{equ:planar:counter}
            \begin{aligned}
        \dot{x}_1 =& -\varepsilon x_1^3+x_2\\
        \dot{x}_2 =& -x_1^3 - x_2^3,
    \end{aligned}
    \end{equation}
 for some $\varepsilon> 0$. One can use the Lyapunov function $\mathbb{R}^2 \ni (x_1,x_2)\mapsto {V}(x_1,x_2):=x_1^4/4+x_2^2/2$ to check that $x^*=0$ is indeed GAS. However, one can show that $V$ is not Morse, and in fact, no Morse Lyapunov function exists for \eqref{equ:planar:counter} \cite[Prop. 3.5]{ref:jongeneel2026global}. It seems reasonable that hyperbolicity plays a role here, however, the real problem is having a linearization where the imaginary part cannot be diagonalized \cite[Cor. 3.9]{ref:jongeneel2026global}. We generalize this observation below, that is, not only does a locally quadratic Lyapunov function fail for \eqref{equ:planar:counter}, even a Lyapunov function that is sandwiched between positively homogeneous polynomials, in $\|x\|$, of equal degree (\textit{e.g.}, quadratic forms), fails to exist \textit{cf}. the quadratic bounds in \cite{ref:zkolter2019learning}.  

 \begin{proposition}[An obstruction to local homogeneity]
 \label{prop:obstruction}
     Let $f\in \mathcal{S}^r_0(\mathbb{R}^n)$, for some $r\in \mathbb{N}_{\geq 1}\cup \{\infty\}$ and $n\in \mathbb{N}_{\geq 1}$. Suppose that $\mathrm{D}f(0)$ has eigenvalues on the imaginary axis that are not semisimple. Suppose that we pick some $C^0$ function $V:\mathbb{R}^n\to \mathbb{R}_{\geq 0}$ such that
     \begin{equation}
     \label{equ:sandwich:V}
         a_1 \| x\|^p \leq V(x) \leq a_2 \|x\|^p \quad \forall x\in U, 
     \end{equation}
     for some $0<a_1\leq a_2$, $p\in \mathbb{N}_{\geq 1}$ and $U$ being some neighbourhood of $0$. Then, $V$ cannot be a Lyapunov function for $x^*=0$ under $f$. 
 \end{proposition}
\begin{proof}
    Since $\mathrm{D}f(0)$ has an eigenvalue that is not semisimple, there is some $x_0\in U$ such that $\lim_{t\to +\infty}\|e^{t\mathrm{D}f(0)}x_0\|= +\infty$. This implies in particular that there is some $T>0$ such that $a_1 \| e^{T\mathrm{D}f(0)}x_0 \|^p > a_2 \| x_0\|^p$. However, then for a sufficiently small $\delta>0$ (\textit{i.e.}, to stay within $U$) one has
    \begin{align*}
        V(\phi^T(\delta x_0)) \geq a_1 \| \phi^T(\delta x_0)\|^p = a_1 \| \delta e^{T \mathrm{D}f(0)} x_0 + r(\delta)\|^p,
    \end{align*}
    where $\delta \mapsto r(\delta) $ is a remainder term such that $\lim_{\delta \to 0^+} r(\delta)/\delta = 0$. 
    Continuing, after a further decrease in $\delta>0$, if necessary, one obtains:  
    \begin{align*}
         a_1 \delta^p \| e^{T \mathrm{D}f(0)} x_0 + r(\delta)/\delta\|^p > a_2 \delta^p \| x_0\|^p \geq V(\delta x_0). 
    \end{align*}
    Hence, $V$ is increasing along a solution and so it cannot be a valid Lyapunov function. 
\end{proof}
Note, one may continue this line of arguments and generalize \eqref{equ:sandwich:V} to a subset of $\mathcal{K}_{\infty}$ functions. We leave this, and in fact a complete characterization, for future work. 

A typical non-hyperbolic example that \textit{does} admit a Morse Lyapunov function is of the form
  \begin{equation}
        \label{equ:planar:nonhyperbolic:with:Morse}
            \begin{aligned}
        \dot{x}_1 =& x_2 - x_1(x_1^2+x_2^2),\\
        \dot{x}_2 =& -x_1-x_2(x_1^2+x_2^2),
    \end{aligned}
    \end{equation}
    \textit{e.g.}, see~\cite[\S 7.6]{Sastry} and \cite[Ex. 3.7]{ref:jongeneel2026global}. One can check that $(x_1,x_2)\mapsto V(x_1,x_2):=x_1^2+x_2^2$ is a Morse Lyapunov function indeed. In particular, $\dot{V}(x)=-2\|x\|_2^4$. 
 
\section{Morse Lyapunov functions}
Inspired by genericity results (\textit{e.g.}, see \cite{ref:peixoto1967approximation}, \cite[p. 147]{ref:hirsch1976differential} or \cite[p. 153]{palis1982geometric}), one might expect that Morse Lyapunov functions still provide for a sufficiently rich class of functions. To that end, let $\mathcal{HS}_0^r(\mathbb{R}^n)\subset \mathcal{S}_0^r(\mathbb{R}^n)$ be the subset of vector fields such that $0\in \mathbb{R}^n$ is a hyperbolic equilibrium point (\textit{i.e.}, $\Re\{\lambda_i(\mathrm{D}f(0))\}\neq 0$ for $i\in [n]$). In our companion paper we have provided the following converse result, confirming part of what one might expect.  

\begin{theorem}[On the existence of Morse Lyapunov functions {\cite[Thm. 3.1]{ref:jongeneel2026global}}]
\label{thm:Morse:exists:hyperbolic}
  Let $f\in \mathcal{HS}^r_0(\mathbb{R}^n)$, for some $r\in \mathbb{N}_{\geq 1}\cup \{\infty\}$ and $n\in \mathbb{N}_{\geq 1}$, then, there is a smooth Morse Lyapunov function certifying that $0$ is GAS. 
\end{theorem}

A generalization of the Morse Lemma \cite[Lem. 2.2]{ref:milnor1963morse} makes Theorem \ref{thm:Morse:exists:hyperbolic} interesting, that is, through \cite{ref:boumal2026smooth} it follows that Morse Lyapunov functions are elements of 
\begin{equation*}
    \mathscr{V} := \big\{ V\in C^{\infty}(\mathbb{R}^n;\mathbb{R}_{\geq 0})\,|\, \exists\, \varphi\in \mathrm{Diff}^{\infty}(\mathbb{R}^n;\mathbb{R}^n):V(\varphi^{-1}(y))=\inf_{y'\in \mathbb{R}^n} V(y')+\|y\|_2^2\big\}, 
\end{equation*}
where $\mathrm{Diff}^r(\mathbb{R}^n;\mathbb{R}^n)$ denotes the group of $C^r$-smooth diffeomorphisms from $\mathbb{R}^n$ to itself ($C^r$-smooth bijective maps with a $C^r$-smooth inverse). 
In fact, one may pick $\varphi \in \mathrm{Diff}^{\infty,+}(\mathbb{R}^n;\mathbb{R}^n)$ thanks to symmetry in $y\mapsto \| y\|_2^2$ (\textit{e.g.}, see \cite[Prop. 2.1]{ref:jongeneel2026global}), that is, the diffeomorphism can be selected to be orientation preserving (like time-$t$ maps $x\mapsto \phi^t (x)$, for $\phi$ a flow). One may spot that results of this form are closely related to the study of linearizing nonlinear dynamical systems. In general, a diffeomorphism does not suffice, one needs to sacrifice regularity at the equilibrium, see \cite{ref:grune1999asymptotic,ref:JongeneelSCL26,ref:kvalheim2025global}. 

\begin{remark}[Regions of attraction that are not $\mathbb{R}^n$]
\label{rem:ROA}
\begin{upshape}
We have focussed on global asymptotic stability on $\mathbb{R}^n$ to keep the exposition simple. However, consider some smooth $n$-dimensional manifold $M$ and let $p^*\in M$ be asymptotically stable under some smooth vector field $f$ on $M$, with $B^*$ being the corresponding region (basin) of attraction. A topological fact is that $B^*$ is diffeomorphic to $\mathbb{R}^n$~\cite[Thm. 2.2]{ref:wilson1967structure}. Hence, effectively anything we discuss carries over.  
\end{upshape}
\end{remark}

Then, inspired by \cite{ref:wilson1967structure}, \cite{ref:cheng2024learning} proposed to use a function $g_{\theta}$, parametrized via \textsf{BiLipNet} (or \textsf{PLNet} for that matter) \textit{e.g.}, see \cite{ref:pmlr-v235-wang24p}, and study Lyapunov functions of the form $x\mapsto V_{\theta}(x)=\|g_{\theta}(x)\|_2^2$, which we understand as a predecessor of using $\mathscr{V}$ (as global smooth P\L{} (Lyapunov) functions, \textit{i.e.}, $\|\nabla V_{\theta}(x)\|^2_2\geq \mu V_{\theta}(x)$, for some $\mu>0$ and all $x\in \mathbb{R}^n$, are a subset of $\mathscr{V}$ \cite[Thm. 1.1]{ref:boumal2026smooth}). They identified that such a parametrization directly aids in enforcing structure relevant for Lyapunov functions. Since $g_{\theta}$ is bi-Lipschitz, the map is invertible, and in contrast to other work, $V_{\theta}$ is not enforced to be convex. Convexity may be restrictive \citep{ref:JongeneelSchwan2024TAC}. Then, motivated by finding neural architectures that guarantee $\nabla V_{\theta}(x) \neq 0$ for all $x\neq x^*$, \cite{ref:zhong2026polarnet} proposed \textsf{PolarNet} in a recent preprint. They effectively recover a subset of the results from \cite{ref:boumal2026smooth}, \textit{cf}. \cite[Thm. 1]{ref:zhong2026polarnet}. They work with $x\mapsto V_{\theta}(x)=\|\psi_{\theta}(x)\|_2^2$ for an appropriate map $\psi_{\theta}:\mathbb{R}^n\to \mathbb{R}^n$ and claim that if $ \psi_{\theta}$ is an \textit{affine coupling flow architecture}, then any $C^2$ diffeomorphism can be approximated, up to arbitrary precision. The authors appeal to \cite{ref:teshima2020coupling} to substantiate this claim. However, that work provides a universal approximation property with respect to the $L_p$ or $\mathrm{sup}$-norm; and these norms are not of use in the context of Lyapunov inequalities. Then, somewhat paraphrasing, \cite{ref:zhong2026polarnet} write that asymptotic stability implies that a Lyapunov function exists that can be accurately modelled using \textsf{PolarNet}. Below we formalize and sharpen this claim.

\section{Topologies to certify stability}
Merely universal approximation in \textit{some} norm does not readily suffice in the context of Lyapunov functions. Suppose that $V$ is an unknown Lyapunov function with respect to some dynamical system \eqref{equ:f}, then any $\widetilde{V}$ that is close, in the sup-norm, to $V$ is by no means a valid Lyapunov function. To get a grip on regularity, define the $C^k$-norm of a $C^k$-smooth map $h:\mathbb{R}^d\to \mathbb{R}^n$, on some (bounded) domain $\Omega\subseteq \mathbb{R}^d$, for simplicity, through
\begin{align*}
    \| h\|_{C^k(\Omega)} := \sup_{i\in \mathbb{N}_{[0,k]}}\sup_{x\in \Omega}\| \mathrm{D}^i h(x)\|. 
\end{align*}

The standard example to have in mind is as follows, given a function $h\in C^1(\mathbb{R};\mathbb{R})$, then a perturbation of the form $\mathbb{R}\ni x\mapsto e(x):=\epsilon \sin(x/\epsilon^2)$ results in $\|h+e\|_{C^0([0,1])}\leq \|h\|_{C^0([0,1])}+\epsilon$. However, $\|h+e\|_{C^1([0,1])}\geq h'(0)+1/\epsilon$. In the context of Lyapunov functions, merely $C^0$ control is futile, \textit{e.g.}, $\| V - V_{\theta}\|_{C^0(K)}<\epsilon$, for a small $\epsilon>0$, does not imply that $L_f V_{\theta}<0$, despite $L_f V<0$.

Without going into measure theory, the $L_p$-norm of such a function $h$ is defined through $\|h\|_{L_p(\Omega)} := (\int_{\Omega} \|h\|^p \mathrm{d}\mu)^{1/p}$ for $p\in \mathbb{N}_{\geq 1}$. The sup-norm, on the other hand, is simply defined as $\|h\|_{\Omega,\infty}:=\sup_{x\in \Omega} \| h(x)\|$. It is now important to note that both the $L_p$-norm and the sup-norm do not control the $C^k$-norm for $k>0$. Hence, they are of little use in the context of learning Lyapunov functions.   

This is known, \cite{ref:harting2026locally} consider $C^1$-control as well, although not directly by construction, but by appealing to standard universal approximation results.  

\section{Universal approximation properties of invertible neural networks}
\label{sec:universal:approximation}
\cite{ref:zhong2026polarnet} appeal to \cite{ref:teshima2020coupling} and claim that $L_p$ or sup-norm control suffices to learn Lyapunov functions. We will appeal to more recent results by \cite{ref:ishikawa2023universal} regarding invertible neural networks (INNs) and control of a Sobolev norm. To follow \cite{ref:ishikawa2023universal}, we need to work with $n\geq 2$. We return to this in Remark~\ref{rem:n:1}. 

Let $\mathcal{G}\subseteq \mathrm{Lip}^r(\mathbb{R}^n;\mathbb{R}^n)$ (\textit{i.e.}, some $C^r$-smooth globally Lipschitz maps) and let $\phi_g$ be the flow corresponding to $\dot{y}=g(y)$ for any $g\in \mathcal{G}$. Then define the set of diffeomorphisms $\Phi(\mathcal{G}):=\{ x\mapsto \phi_g^1(x):g\in \mathcal{G}\}$. Now we may define a set of \textbf{\textit{invertible neural networks}} (INNs) through
\begin{align*}
    \mathrm{INN}_{\Phi(\mathcal{G})}:=\{W_1\circ \varphi_1\circ \cdots \circ W_k\circ \varphi_k: k\in \mathbb{N}_{\geq 1}, \varphi_i \in \Phi(\mathcal{G}),W_i\in \mathsf{Aff}(\mathbb{R}^n)\},
\end{align*}
where $\mathsf{Aff}(\mathbb{R}^n)$ denotes the affine group, \textit{i.e.}, $x\mapsto W_i(x)=A_ix+b_i$ for some $A_i\in \mathsf{GL}(\mathbb{R}^n)$ and $b_i\in \mathbb{R}^n$. Importantly, $\mathrm{INN}_{\Phi(\mathcal{G})}\subseteq \mathrm{Diff}^{r}(\mathbb{R}^n;\mathbb{R}^n)$. Moreover, by composing with some appropriate $W_0\in \mathsf{Aff}(\mathbb{R}^n)$, one can always force the resulting diffeomorphism to be such that $0$ maps to itself.  

Then, we say that a set of maps $\widehat{\mathcal{G}}\subseteq C^r(\mathbb{R}^n;\mathbb{R}^n)$ is a $W^{r,\infty}$-universal approximator for $\mathcal{G}\subseteq C^r(\mathbb{R}^n;\mathbb{R}^n)$ if for any compact set $K\subseteq \mathbb{R}^n$, any $\epsilon>0$ and any $g\in \mathcal{G}$ there is a $\widehat{g}\in \widehat{\mathcal{G}}$ such that
\begin{align*}
    \|g-\widehat{g}\|_{K,r,\infty} := \sum_{|\alpha|\leq r} \sup_{x\in K} \| \mathrm{D}^{\alpha} g(x)-\mathrm{D}^{\alpha} \widehat{g}(x)\| < \epsilon.
\end{align*}

The next theorem characterizes $W^{r,\infty}$-universality of \textit{neural ODEs} (NODEs), that is, using NODEs, one has control over a Sobolev norm. 
\begin{theorem}{{\cite[Thm. 44]{ref:ishikawa2023universal}}}
\label{thm:NODE:universality}
    Let $r\in \mathbb{N}_{\geq 0}$ and suppose that $\widehat{\mathcal{G}}\subseteq \mathrm{Lip}^r(\mathbb{R}^n;\mathbb{R}^n)$ is a $W^{r,\infty}$-universal approximator for $\mathrm{Lip}^r(\mathbb{R}^n;\mathbb{R}^n)$, then $\mathrm{INN}_{\Phi(\widehat{\mathcal{G}})}$ is a $W^{r,\infty}$-universal approximator for $\mathrm{Diff}^{\mathrm{max}\{r,1\}}(\mathbb{R}^n;\mathbb{R}^n)$. 
\end{theorem}

Concretely, what Theorem~\ref{thm:NODE:universality} means is that if we fix some $r\geq 0$, then for some appropriate set of globally Lipschitz vector fields $\widehat{\mathcal{G}}$, one may hand us any diffeomorphism $\varphi\in \mathrm{Diff}^{\mathrm{max}\{r,1\}}(\mathbb{R}^n;\mathbb{R}^n)$, a compact set $K\subseteq \mathbb{R}^n$ and some $\epsilon>0$, then we know that there is some $k\in \mathbb{N}_{\geq 1}$ such that
\begin{align}
\label{equ:Sobolev:eps:bound}
    \| \varphi - W_1\circ \varphi_1\circ \cdots \circ W_k\circ \varphi_k \|_{K,r,\infty} < \epsilon 
\end{align}
for $\varphi_1,\dots,\varphi_k\in \Phi(\widehat{\mathcal{G}})$ and $W_1,\dots,W_k\in \mathsf{Aff}(\mathbb{R}^n)$. Appropriate choices of $\widehat{\mathcal{G}}$ are discussed in \cite[p. 19]{ref:ishikawa2023universal}, \textit{e.g.}, \textsf{MLPs} with finite weights and $\mathsf{tanh}$ activation functions (since we typically consider $r>0$).  

As mentioned before, one can always find an $\mathrm{INN}_{\Phi(\widehat{\mathcal{G}})}$-generated diffeomorphism $\varphi_{\theta}$ that maps $0$ to $0$. Regarding control, let $K$ be a compact neighbourhood of $0$ and let $\varphi_{\theta}:=W_1\circ \varphi_1\circ \cdots \circ W_k\circ \varphi_k$ be such that \eqref{equ:Sobolev:eps:bound} holds for $\varphi$ generating some $V\in \mathscr{V}$. Let $x\mapsto W_0(x):=x-b_0$ for $b_0:=\varphi_{\theta}(0)$, then $W_0\circ \varphi_{\theta}$ (\textit{e.g.}, adjust $W_1$) is a diffeomorphism and such that  
 \begin{align*}
    \| \varphi - W_0\circ \varphi_{\theta}\|_{K,r,\infty} = \| \varphi - W_1\circ \varphi_1\circ \cdots \circ W_k\circ \varphi_k +b_0 \|_{K,r,\infty} < \epsilon + \|b_0\|_{K,r,\infty}<2\epsilon 
\end{align*}
as \eqref{equ:Sobolev:eps:bound} implies that $\|\varphi(0)-\varphi_{\theta}(0)\|<\epsilon$ and $\varphi(0)=0$ as $\varphi$ generates some $V\in \mathscr{V}$. 
\section{Approximation of stability certificates}
Combining Theorem~\ref{thm:Morse:exists:hyperbolic} and Theorem~\ref{thm:NODE:universality}, we can detail approximation quality of NODE-based INNs (\textit{i.e.}, $\mathrm{INN}_{\Phi(\mathcal{G})}$) to construct Lyapunov functions. 
Regarding notation, we recall that $\|g\|_{K,\infty}:=\|g\|_{K,0,\infty}=\sup_{x\in K} \| g(x)\|$. Regarding the power of INNs, one should observe that below we effectively only need to be concerned with $L_f V_{\theta}<0$, other properties are direct, globally. 

\begin{theorem}[Universal approximation of Morse Lyapunov functions for hyperbolic equilibria]
\label{thm:main}
    Let $f\in \mathcal{HS}^r_0(\mathbb{R}^n)$, for some $r\in \mathbb{N}_{\geq 1}\cup \{\infty\}$ and $n\in \mathbb{N}_{\geq 2}$, then
    \begin{enumerate}[(i)]
        \item there is a Lyapunov function $V\in \mathscr{V}$ (\textit{i.e.}, $V=\|\varphi\|_2^2$ for $\varphi\in \mathrm{Diff}^{\infty}(\mathbb{R}^n;\mathbb{R}^n)$) such that $PA+A^{\mathsf{T}}P\prec 0$ for $A:=\mathrm{D}f(0)$, $P:=\nabla^2 V(0)$, certifying that $x^*=0$ is GAS;
        \item    for any compact neighbourhood $K$ of $0$ there is an $\epsilon_K>0$ such that if some diffeomorphism $\varphi_{\theta}\in \mathrm{Diff}^1(\mathbb{R}^n;\mathbb{R}^n)$ satisfies $\|\varphi_{\theta}-\varphi\|_{K,1,\infty}<\epsilon_K$, then $x\mapsto V_{\theta}(x):=\|\varphi_{\theta}(x)-\varphi_{\theta}(0)\|_2^2$ is a Lyapunov function satisfying $L_f V_{\theta}<0$ on $K\setminus \{x^*\}$, thus certifying global asymptotic stability on $\{x\in K:V_{\theta}(x)\leq \min_{x\in \partial K} V_{\theta}(x)\}$; and
        \item such a $\varphi_{\theta}$ can be generated through an invertible neural network $\mathrm{INN}_{\Phi(\widehat{\mathcal{G}})}$, where $\widehat{\mathcal{G}}\subseteq \mathrm{Lip}^1(\mathbb{R}^n;\mathbb{R}^n)$ is a $W^{1,\infty}$-universal approximator of  $\mathrm{Lip}^1(\mathbb{R}^n;\mathbb{R}^n)$.  
    \end{enumerate}
\end{theorem}
\begin{proof}
    \textit{(i)} Since $f\in \mathcal{HS}^r_0(\mathbb{R}^n)$, by Theorem~\ref{thm:Morse:exists:hyperbolic} there is a Lyapunov function $V\in \mathscr{V}$ certifying that $x^*=0$ is GAS, \textit{i.e.}, there is a $\varphi \in \mathrm{Diff}^{\infty,+}(\mathbb{R}^n;\mathbb{R}^n)$ such that $x\mapsto V(x) = \| \varphi(x)\|_2^2$. The fact that $V$ can be chosen to satisfy a Lyapunov equation (or inequality $PA+A^{\mathsf{T}}P\prec 0$ for that matter), locally, follows from the proof of \cite[Thm. 3.1]{ref:jongeneel2026global}.

    \textit{(iii)} By Theorem~\ref{thm:NODE:universality}, we can find a $\varphi_{\theta}\in \mathrm{Diff}^1(\mathbb{R}^n;\mathbb{R}^n)$, generated by some invertible neural network $\mathrm{INN}_{\Phi(\widehat{\mathcal{G}})}$, such that $\| \varphi - \varphi_{\theta}\|_{K,1,\infty}<\epsilon_K$. 

    \textit{(ii)} We claim that there is a choice of $\epsilon_K$ such that $x\mapsto \| \varphi_{\theta}(x)\|^2_2=:V_{\theta}(x)$ is a Lyapunov function on $K$. As discussed in Section~\ref{sec:universal:approximation}, if needed, one can redefine the map $\varphi_{\theta}$ to be $x\mapsto \varphi_{\theta}(x)-\varphi_{\theta}(0)$, \textit{i.e.}, via composition with some $W_0\in \mathsf{Aff}(\mathbb{R}^n)$. Hence $V_{\theta}(x)=0\iff x=0$. We remark, with some abuse of notation, that if $\| \varphi - \varphi_{\theta}\|_{K,1,\infty}<\tfrac{1}{2}\epsilon_K$, then the shifted map satisfies $\| \varphi - \varphi_{\theta}\|_{K,1,\infty}<\epsilon_K$.

    Secondly, we see that $V_{\theta}(x)>0$, for all $x\in K\setminus \{0\}$ and in fact, all sublevel sets $\{ x\in \mathbb{R}^n:V_{\theta}(x)\leq c\}$, for $c>0$, are compact as those sets are diffeomorphisms of the compact balls $\overline{B_{\sqrt{c}}(0)}=\{x\in \mathbb{R}^n: \|x\|^2_2 \leq c\}$.

    The remaining property to check is that $L_f V_{\theta}<0$ on $K\setminus \{0\}$. We know that $L_fV<0$ on $K\setminus \{0\}$, plus that $x\mapsto L_f V_{\theta}(x)= L_f \|\varphi_{\theta}(x)\|_2^2= 2\langle \varphi_{\theta}(x), \mathrm{D}\varphi_{\theta}(x) f(x) \rangle$. Consider now some compact ball $\overline{B}_{\varrho}(0)\subseteq \mathrm{int}(K)$, for some $\varrho>0$ that we select later, and study $L_f V_{\theta}$ on the compact set $W:=K\setminus B_{\varrho}(0)$. Let $-2\beta := \max_{x\in W}L_f V(x)<0$. Then, as $L_f V_{\theta}=L_f V +(L_fV_{\theta}-L_fV)$. We have that $L_f V_{\theta}<0$ on $W$ if $(L_fV_{\theta}-L_fV)\leq \beta$ on $W$. Now, as $\mathrm{D}\varphi^{\mathsf{T}}\varphi-\mathrm{D}\varphi_{\theta}^{\mathsf{T}}\varphi_{\theta} = (\mathrm{D}\varphi-\mathrm{D}\varphi_{\theta})^{\mathsf{T}}\varphi+\mathrm{D}\varphi_{\theta}^{\mathsf{T}}(\varphi-\varphi_{\theta})$ and $\| \mathrm{D}\varphi_{\theta}\|_{W,\infty}\leq \| \mathrm{D}\varphi\|_{W,\infty}+\| \mathrm{D}\varphi_{\theta}-\mathrm{D}\varphi\|_{W,\infty}$ see that 
    \begin{align*}
        |L_fV_{\theta}-L_fV| = |L_f(V-V_{\theta})| &\leq \|f\|_{W,\infty}\| \mathrm{D}V-\mathrm{D}V_{\theta}\|_{W,\infty}\\
        &\leq 2\|f\|_{W,\infty}\| \mathrm{D}\varphi^{\mathsf{T}}\varphi-\mathrm{D}\varphi_{\theta}^{\mathsf{T}}\varphi_{\theta}\|_{W,\infty}\\
        &\leq 2\|f\|_{W,\infty}\| \mathrm{D}\varphi-\mathrm{D}\varphi_{\theta}\|_{W,\infty}\|\varphi\|_{W,\infty}\\
        &+ 2\|f\|_{W,\infty}\| \varphi-\varphi_{\theta}\|_{W,\infty}\left( \|\mathrm{D}\varphi\|_{W,\infty}+\| \mathrm{D}\varphi_{\theta}-\mathrm{D}\varphi\|_{W,\infty}  \right).
    \end{align*}
Let $\kappa_1:=\|f\|_{W,\infty}$, $\kappa_2:=\|\varphi\|_{W,\infty}$ and $\kappa_3:=\|\mathrm{D}\varphi\|_{W,\infty}$. If $\|\varphi-\varphi_{\theta}\|_{W,1,\infty}\leq \epsilon$ holds for
\begin{align*}
    0 < \epsilon \leq \frac{-(\kappa_2+\kappa_3) + \sqrt{(\kappa_2+\kappa_3)^2 +2\beta/\kappa_1}}{2},
\end{align*}
then  $|L_fV_{\theta}-L_fV|\leq \beta$. Hence, $\max_{x\in W}L_f V_{\theta}\leq -\beta < 0$. 

Now we discuss what happens within $\overline{B_{\varrho}(0)}$, better yet, how to select $\varrho>0$. Observe that $P=\nabla^2 V(x^*)=2 \mathrm{D}\varphi(0)^{\mathsf{T}}\mathrm{D}\varphi(0)$, and this $P$ satisfies $P\mathrm{D}f(0)+\mathrm{D}f(0)^{\mathsf{T}}P\prec0$. Following a standard argument (\textit{e.g.}, see the proof of \cite[Thm. 4.7]{ref:khalil2002nonlinear}), it can be shown that there is a $\varrho>0$ and $\kappa_{\varrho}>0$ such that $L_f V(x)<-2\kappa_{\varrho} \|x\|_2^2$ on $B_{\varrho}(0)$. Indeed, the reason that such a $V$ exists follows since $f\in \mathcal{HS}^r_0(\mathbb{R}^n)$, see \cite[Thm. 3.1]{ref:jongeneel2026global}.

    Then, as above, a useful decomposition is as follows:
    \begin{align*}
        \tfrac{1}{2}L_f (V_{\theta}- V) = \langle \varphi_{\theta}-\varphi, \mathrm{D}\varphi f\rangle + \langle \varphi, \mathrm{D}(\varphi_{\theta}-\varphi) f\rangle + \langle \varphi_{\theta}-\varphi, \mathrm{D}(\varphi_{\theta}-\varphi) f\rangle. 
    \end{align*}
    Akin to $\kappa_{\varrho}$, one can set $\kappa_f := \sup_{x\in B_{\varrho}(0)\setminus \{0\}} \|f(x)\|_2/\|x\|_2$, where one may shrink $\varrho>0$ if needed. In addition, set $\kappa_4:=\|\mathrm{D}\varphi\|_{\overline{B_{\varrho}(0)},\infty}$. Now, if the Sobolev bound $\|\varphi-\varphi_{\theta}\|_{\overline{B_{\varrho}(0)},1,\infty}\leq \epsilon'$ holds, one obtains, via the mean value theorem and convexity of $B_{\varrho}(0)$ (\textit{e.g.}, $\|\varphi_{\theta}(x)-\varphi(x)\|_2 \leq \epsilon' \|x\|_2$), that $|L_f(V_{\theta}-V)(x)|\leq 2\kappa_f \epsilon' (2\kappa_4 + \epsilon') \|x\|_2^2$. Hence, if $\epsilon'>0$ is such that
\begin{align*}
    0 < \epsilon' \leq \frac{-2\kappa_4 + \sqrt{4\kappa_4^2 +2\kappa_{\varrho}/\kappa_f}}{2},
\end{align*}
then  $|L_fV_{\theta}(x)-L_fV(x)|\leq \kappa_{\varrho} \|x\|_2^2 $. In its turn, this implies that $L_f V_{\theta}(x)<-\kappa_{\varrho} \|x\|_2^2$ on $B_{\varrho}(0)$. 

At last, we conclude by selecting $\epsilon_K = \tfrac{1}{2}\min\{\epsilon,\epsilon'\}$. 
\end{proof}

Theorem~\ref{thm:main} can only handle $n\geq 2$ (as item \textit{(iii)} relies on \cite{ref:ishikawa2023universal}), however, we mention the following fact. As such, the case $n=1$ is irrelevant from the learning perspective. 

\begin{remark}[The case $n=1$]
\label{rem:n:1}
\begin{upshape}
    Let $f\in \mathcal{S}^r_0(\mathbb{R}^n)$,  for some $r\in \mathbb{N}_{\geq 1}\cup \{\infty\}$ and suppose that $n=1$. In this case, the trivial Morse Lyapunov function, \textit{i.e.}, $x\mapsto V(x)=x^2$, always works as $xf(x)<0$ on $\mathbb{R}\setminus \{0\}$. See~\cite{ref:sontag1980remarks} for further remarks. 
\end{upshape}
\end{remark}

The proof of Theorem~\ref{thm:main} reveals that as long as asymptotic stability can be certified using a Morse Lyapunov function, one can still characterize a form of universality. However, on an annulus/shell. 

\begin{corollary}[Universal approximation of Morse Lyapunov functions]
\label{cor:main}
    Let $f\in \mathcal{S}^r_0(\mathbb{R}^n)$, for some $r\in \mathbb{N}_{\geq 1}\cup \{\infty\}$ and $n\in \mathbb{N}_{\geq 2}$. Suppose, there is a Lyapunov function $V\in \mathscr{V}$ (\textit{i.e.}, $V=\|\varphi\|_2^2$ for $\varphi\in \mathrm{Diff}^{\infty}(\mathbb{R}^n;\mathbb{R}^n)$), certifying that $x^*=0$ is GAS. Then, 
    \begin{enumerate}[(i)]
        \item    for any compact neighbourhood $K$ of $0$ and radius $\varrho>0$ such that $\overline{B_{\varrho}(0)}\subseteq \mathrm{int}(K)$, there is an $\epsilon_{K,\varrho}>0$ such that if some diffeomorphism $\varphi_{\theta}\in \mathrm{Diff}^1(\mathbb{R}^n;\mathbb{R}^n)$ satisfies $\|\varphi_{\theta}-\varphi\|_{K,1,\infty}<\epsilon_{K,\varrho}$, then $x\mapsto V_{\theta}(x):=\|\varphi_{\theta}(x)-\varphi_{\theta}(0)\|_2^2$ is a Lyapunov function satisfying $L_f V_{\theta}<0$ on $K\setminus B_{\varrho}(0)$; and
        \item such a $\varphi_{\theta}$ can be generated through an invertible neural network $\mathrm{INN}_{\Phi(\widehat{\mathcal{G}})}$, where $\widehat{\mathcal{G}}\subseteq \mathrm{Lip}^1(\mathbb{R}^n;\mathbb{R}^n)$ is a $W^{1,\infty}$-universal approximator of  $\mathrm{Lip}^1(\mathbb{R}^n;\mathbb{R}^n)$.  
    \end{enumerate}
\end{corollary}
We interpret Corollary~\ref{cor:main} as showing that one can learn a Morse Lyapunov function that certifies \textit{practical} asymptotic stability, \textit{cf}. \cite[\S 25]{ref:lefschetzlasalle1961}. Arguably the most powerful ramification follows when $f\in \mathcal{HS}^r_0(\mathbb{R}^n)$, as in that case, there is always a Lyapunov function $V\in \mathscr{V}$. To include $B_{\varrho}(0)$, one needs to handle that $L_f V_{\theta}$ vanishes at $x^*=0$. To that end, Theorem~\ref{thm:main} \textit{(i)} requires more control of $V$, locally. Such a $V$ always exists \cite[Thm. 3.1]{ref:jongeneel2026global}, but we do not enforce this structure through $\mathrm{INN}_{\Phi(\widehat{\mathcal{G}})}$. Hence, Corollary~\ref{cor:main} provides at the moment the most practically relevant guarantees when using vanilla NODEs and $\mathrm{INN}_{\Phi(\widehat{\mathcal{G}})}$. We are not aware of similar results for other contemporary architectures. However, the following remark illustrates that the above extends to P\L{}-functions. Exploiting this is future work. 

\begin{remark}[Equivalence of Morse and P\L{} Lyapunov functions on compact sets]
\label{rem:PL}
    \begin{upshape}
        It turns out that on compact sets, using smooth  $\mathrm{P}\L{}$ Lyapunov functions is equivalent to using Morse Lyapunov functions. First, $C^2$-smooth $\mathrm{P}\L{}$-functions, with a unique critical point, are Morse by \cite[Cor. 2.17]{ref:rebjock2025fast} and \cite[Thm. 1.1]{ref:boumal2026smooth}. The other direction also works, on compact sets. The local argument follows from \cite[Cor. 2.11]{ref:rebjock2025fast} and continuity allows for an extension to arbitrary compact neighbourhoods. To be somewhat complete, suppose $V$ is Morse. Then, around $x^*=0$ we have $V(x)=\frac{1}{2}\langle x, Px\rangle + o(\|x\|_2^2)$ for $P\succ0$ and $\nabla V(x) = Px + o(\| x\|_2)$, but then there is a $\varrho>0$ such that on $B_{\varrho}(0)$ we have $V(x)\leq 2 \lambda_{\mathrm{max}}(P)\|x\|_2^2$ and $\|\nabla V(x)\|_2\geq \tfrac{1}{2}\lambda_{\mathrm{min}}(P)\|x\|_2$. Putting it together yields $\| \nabla V(x)\|_2^2 \geq \mu_{\varrho} V(x)$ for $\mu_{\varrho}:=(1/8) \lambda_{\mathrm{min}}(P)^2/\lambda_{\mathrm{max}}(P)>0$. Now, if we work on a compact neighbourhood $K$ of $0$, we can define $\mu_1 := \max_{x\in K\setminus B_{\varrho}(0)} V(x)>0$ and $\mu_2:=\min_{x\in K\setminus B_{\varrho}(0)} \|\nabla V(x)\|_2^2>0$ to obtain $\|\nabla V(x)\|_2^2 \geq (\mu_2/\mu_1) V(x)$ on $K\setminus B_{\varrho}(0)$. We can conclude by setting $\mu:=\min\{ \mu_{\varrho},\mu_2/\mu_1\}>0$. 
    \end{upshape}
\end{remark}

\section{Numerical experiments}

\begin{figure}
    \centering
    \includegraphics[width=0.95\linewidth]{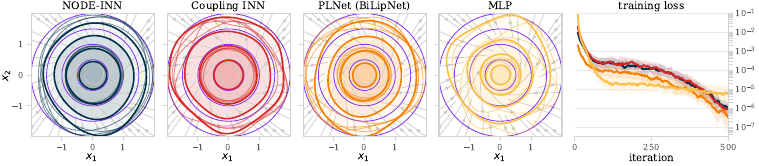}
    \caption{Experiment A. \textsf{NODE-INN} generates $\varphi_{\theta}$ via $\mathrm{INN}_{\Phi(\mathcal{G})}$ with $3$ layers and $\mathsf{tanh}$ activation functions generating Lipschitz vector fields. Numerical integration is done via $\mathsf{RK4}$. \textsf{Coupling-INN} is effectively the same but uses affine coupling flows instead. The \textsf{PLNet} is implemented with $\nu_1=1/2$ and $\nu_2=2$, whereas \textsf{MLP} is a standard $\mathsf{tanh}$ network ($2$ hidden layers with $20$ neurons). We train using $512$ samples, $m=10^{-3}$ and learning rate $10^{-2}$. The purple levelsets correspond to a Lyapunov function for \eqref{equ:planar:nonhyperbolic:with:Morse} and the thinner colour-coded curves are some levelsets of $V_{\theta}$ found after 500 iterations of Adam. We show the outcomes of 100 independent runs. The thick curves visualize the first run. The training loss displays the mean and the 25th-75th percentile.}
    \label{fig:num:1}
\end{figure}

The purpose of this section is to provide an illustration of the above, hence we consider systems that can be visualized (\textit{i.e.}, planar dynamical systems).

We designed the experiments such that all architectures have roughly the same number of parameters. We compare the \textsf{NODE-INN} (\textit{e.g.}, \cite{ref:ishikawa2023universal}
 and this work) against \textsf{Coupling-INN} (\textit{e.g.}, \cite{ref:teshima2020coupling,ref:zhong2026polarnet}), \textsf{PLNet} (via \textsf{BiLipNet}) (\textit{e.g.}, \cite{ref:pmlr-v235-wang24p}) and a standard \textsf{MLP} (\textit{e.g.} using a $\mathsf{tanh}$ activation). The learning is implemented via \textsf{PyTorch}, using a scale-free loss function of the form
 \begin{equation*}
     \theta \mapsto \ell(\theta) := \frac{1}{N}\sum^N_{i=1} \max\left\{0, \frac{L_f V_{\theta}(x_i)}{V_{\theta}(x_i)}+m\right\} 
 \end{equation*}
for some margin $m>0$, samples $\{x_i\}_{i=1}^N$ and $\theta \mapsto V_{\theta}(x):=\| \varphi_{\theta}(x)-\varphi_{\theta}(0)\|_2^2$, where the precise form of $\varphi_{\theta}$ depends on the architecture. We remark that for any $V_{\theta}$ generated through \textsf{NODE-INN}, positive definiteness, properness and the absence of spurious critical points hold by construction. Only the condition $L_f V_{\theta}<0$ needs to generalize beyond the samples. In contrast, for \textsf{MLP}, all defining properties of a Lyapunov function need to generalize.  

\subsection{Experiment A}
We start with learning a Lyapunov function for \eqref{equ:planar:nonhyperbolic:with:Morse}. We present some experimental results in Figure~\ref{fig:num:1}. One observes similar behaviour when comparing \textsf{NODE-INN} and \textsf{PLNet}, although \textsf{PLNet} is faster, both per iteration and with respect to convergence as a function of total iterations.

\subsection{Experiment B}
\label{sec:exp:B}
\begin{figure}
    \centering
    \includegraphics[width=0.95\linewidth]{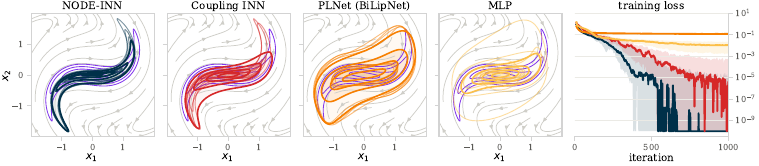}
    \caption{Experiment B. \textsf{NODE-INN}, \textsf{Coupling-INN}, \textsf{PLNet} and \textsf{MLP} generate their respective map $\varphi_{\theta}$ as in Experiment A. The purple levelsets correspond to a Lyapunov function for the nonlinear system from Section~\ref{sec:exp:B}, for $k=6$ and $\beta=1/4$, and the thinner colour-coded curves are some levelsets of $V_{\theta}$ found after 1000 iterations of Adam. We show the outcomes of 30 independent runs. Visualization is as in Figure~\ref{fig:num:1}.}
    \label{fig:num:2}
\end{figure}

\begin{figure}
    \centering
    \includegraphics[width=0.95\linewidth]{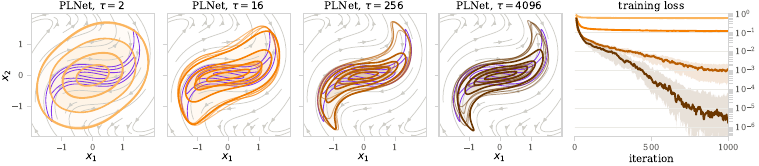}
    \caption{Experiment B, where the distortion $\tau$ with respect to $V_{\theta}$ is adapted.}
    \label{fig:num:3}
\end{figure}

In the second experiment we consider a nonlinear system that is smoothly conjugate to a linear hyperbolic system. In particular, let $\dot{x}=Ax$ and the matrix $J$ be defined through
\begin{equation*}
    A := \begin{pmatrix} -1 & k\\ 0 & -1
    \end{pmatrix},\text{ for some } k\geq 0, \text{ and }  J := \begin{pmatrix}
        0 & -1 \\ 1 & 0
    \end{pmatrix}.
\end{equation*}
Now define $y:=h(x)$ as the time-$1$ map of the flow under $\dot{x}=\beta \|x\|_2^2 J x$, for some $\beta >0$.
Let $x\mapsto V(x)=\|x\|_2^2$ and see that $\langle x,\beta \|x\|_2^2 J x\rangle = 0$. Hence, $h$ is a twist (as the angular velocity is variable). Now, as $A$ is Hurwitz, there is a quadratic form $x\mapsto V_P(x):=\langle x, Px\rangle$ that certifies that $x^*=0$ is GAS under $\dot{x}=Ax$, but then $V_P\circ h^{-1}$ certifies that $y^*=0$ is GAS under $\dot{y}=\mathrm{D}h(h^{-1}(y))Ah^{-1}(y)$. This setting fits well to \textsf{NODE-INN}, and under appropriate choices of $\beta$ and $k$, the nonlinearities can be significant. We present some numerical experiments in Figure~\ref{fig:num:2}. \textsf{NODE-INN} seems to perform well, whereas other architectures appear to resist convergence. However, we recall that with respect to \textsf{BiLipNet}, one employs $(\nu_1,\nu_2)$-Lipschitz maps $\mathcal{F}:\mathbb{R}^n\to \mathbb{R}^n$, that is, $\nu_1 \|x-y\| \leq \| \mathcal{F}(x)-\mathcal{F}(y)\|\leq \nu_2 \| x-y\|$ for all $x,y\in \mathbb{R}^n$. To train the corresponding networks, one sets $\nu_1$ and $\nu_2$, with the \textit{distortion} $\tau:=\nu_2/\nu_1$ (or $(\nu_2/\nu_1)^2$ for $V_{\theta}$) controlling the conditioning, \textit{e.g.}, the condition number of $\nabla^2 V_{\theta}(0)$ is bounded from above by $\tau^2$. The results in Figure~\ref{fig:num:2} are for some particular $\tau$, adapting this parameter may improve performance, as shown in Figure~\ref{fig:num:3}. This is interesting in the light of Remark~\ref{rem:PL}; to work with \textsf{PLNet}, one needs to find an appropriate $\tau$, and this number may be large, \textit{cf}. $\mu_{\varrho}$ from Remark~\ref{rem:PL}. At the same time, this restriction of the search space may aid convergence. A precise characterization is left for future work.

\section{Conclusion}
We have shown that Morse Lyapunov functions are a class of Lyapunov functions with a rich structure that can be exploited both theoretically and algorithmically. Completing the converse theory, also for the other architectures, and providing a rigorous numerical analysis are left for future work.

\subsection*{Generative AI Use}
Generative AI tools (Claude Fable 5.2) aided in comprehending some results in the neural networks literature and those tools were used for code and experimental implementation. The author
reviewed and verified all AI-assisted content and takes responsibility for the
final content of the paper.


\bibliography{l4dc2027-sample}

\end{document}